\documentclass{aims}
\usepackage{txfonts}
\def\typeofarticle{Preprint}

\def\DOI{}
\def\Received{January 2026}
\def\Revised{February 2026}
\def\Accepted{March 2026}
\def\Published{March 2026}

\numberwithin{equation}{section}

\makeatletter
\renewcommand{\@biblabel}[1]{#1\hfill \hspace{-0.2cm}}
\makeatother

\begin{document}

\title{Expanding Flow Shop Tasks Based on Recursive Functions}

\author{%
  Boris Kupriyanov\affil{1}
  Alexander Lazarev\affil{1}
  Aleksandr Roschin\affil{1}
  and
  Frank Werner\affil{2,}\corrauth
}

% \shortauthors is used in copyright information in the end of the paper
\shortauthors{the Author(s)}

\address{%
  \addr{\affilnum{1}}{Trapeznikov Institute of Control Sciences, Russian Academy of Sciences, 117997 Moscow, Russia}
  \addr{\affilnum{2}}{Fakult\"{a}t f\"{u}r Mathematik, Otto-von-Guericke-Universit\"{a}t Magdeburg, PSF 4120. 39016 Magdeburg, Germany}}

% corresponding author
\corraddr{frank.werner@ovgu.de; Tel: +49-391-6752025;\\
 Fax: +49 391 6741171.
}

\begin{abstract}
The paper discusses several extensions of the recursive representation of the flow shop scheduling problem.
It is shown that recursive functions make it possible to describe multiple extensions in a single problem.
The paper considers altogether six extensions.
The examples consider three types of recursive functions: functions associated with the machine, functions that adjust the procession time based on constraints, and functions that control the feasibility of the schedule.
The structure of the superpositions of these functions is presented, and  also
descriptions of several objective functions by recursive functions are presented.
Then the general requirements for a recursive function are formulated and its properties are described.
Finally, a demonstration of the formulation of new problems is provided using examples of simple flow shop extensions and branch and bound optimization.
%\textbf{(200 to 300 words)}
\end{abstract}


\maketitle

\section{Introduction}

The flow shop problem is a classical scheduling problem which has been intensively inivestigated in the literature. A set of $n$ jobs has to be processed on a set of $m$ machines in the same order. Often the permutation flow shop problem (PFSP) is considered, where on all machines the same job order has to be chosen.  This means that one has to select one of the $n!$ possible job permutations such that a certain optimization criterion is satisfied, e.g. the maximum completion time of a job. The PFSP has a number of practical applications in the industry.

This paper presents an original method for describing several extensions of the PFSP.
We describe an approach which is  based on the use of recursive functions.  
The arguments of the recursive function are a pair (job number, machine number), and the calculated value is the completion time of the given job on the given machine.
Subsequently, we consider three types of recursive functions: functions associated with a machine, functions that adjust %åðó 
processing time based on constraints, and functions that control the feasibility of a schedule.
The structure of the superpositions of these functions is described.
In total, this paper discusses six types of PFSP extensions that have not yet been considered as a single problem. 

The consideration of recursive functions in connection with scheduling problems is rather new.
The paper \cite{KLRW} discusses the use of recursive functions to extend the PFSP with the predicate {\em and} and also in connection with job shop problems, where the jobs may have different machine routes. 

The applicability of recursive functions is determined by two of their main features. 
The first and main feature of a recursive function is that it is, by definition, calculated using some algorithm \cite{St,HS}.
This distinguishes it from an analytical function.
The same function can be calculated using several algorithms.
This feature allows, when solving problems in scheduling theory, to localize a part of the algorithmic calculations in the body of a recursive function and make them basic objects for describing and solving certain problems.
The second feature is the ability to define a recursive function through itself and other functions, which allows for complex structural calculations, for example on a precedence graph which, in another setting, are only possible using cyclic calculations.
The applicability of recursive functions to scheduling problems is due to the fact that precedence relations are defined for them both between the machines and between the jobs (with a fixed order of jobs) \cite{KLRW}.

Most of the problems considered in this paper are NP-hard, and exact solutions are possible only for small dimensions.
Among the methods used to find approximate solutions, such as, for example, the ant colony (AC) method \cite{Do} or the simulated annealing (SA) \cite{He} method, for optimizing a certain function, essentially only its computability is required.
In \cite{KR}, the application of one exact and two approximate methods for solving PFSPs is investigated: branch and bound, ant colony, and simulated annealing.
The problems are described by recursive functions and use the predicate {\em and}; this model allows solving planning problems for certain types of assembly production.
The precedence graphs for those problems are generated by a pseudo-random directed acyclic graph generator, which allows for a large number of generated statistics and calculated results. 80 tests were performed for 14, 100 and 200 jobs.
The paper \cite{KR} also discusses a common SDST manufacturing problem with sequence-dependent setup times.

This outline of this paper is as follows:
In Section 2, the description of the  PFSP and the definition of the precedence relation between the jobs and machines is given.
The representation of the PFSP using recursive functions is considered in Section 3.
In Section 4, the description of recursive functions of PFSP extensions and the structure of function superpositions is given.
Section 5 introduces some typical objective functions and in Section 6,
the description of a subset of recursive functions applicable to solving the PFSP and the properties of these functions are given. 
Section 7 presents the description of the branch and bound method for an extended PFSP 
and an example of a problem solution. Finally, Section 8 gives some concluding results.

\section{Description of the flow shop problem class}
 Consider the description of the PFSP  \cite{L2}.
Let\\
${\cal   J}=\{1,2,\dots,n\}$ be the set of jobs;\\    % $Z \; (|Z|=n)$,\\
${\cal M}=\{1,2,\dots,m\}$ be the set of machines arranged sequentially;\\
$P(n,m)$ be the matrix of the job completion times, where $p_{j,k}$ is the
processing time of job $j$ on machine $k$.\\
\begin{figure}[!ht]
	\centering
	\includegraphics [scale=0.5]{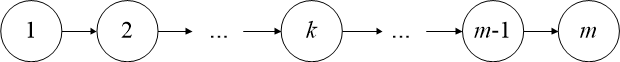}
	\caption{Example of a flow shop problem chain.}
	\label{fig1a}
\end{figure}
The graph of the precedence relations of the machines for the single job (see Figure  \ref{fig1a}) is a simple chain in which the vertices are machines and the arcs define the precedence relationships between them.
With a certain fixed order of jobs, this graph can be ''expanded'', and the precedence relationships between the jobs ${\cal J}\times {\cal J}$ and between the machines ${\cal M}\times {\cal M}$ (see Figure  \ref{fig2}) can be defined.

 \begin{figure}[!ht]
 	\centering
 	\includegraphics [scale=0.6]{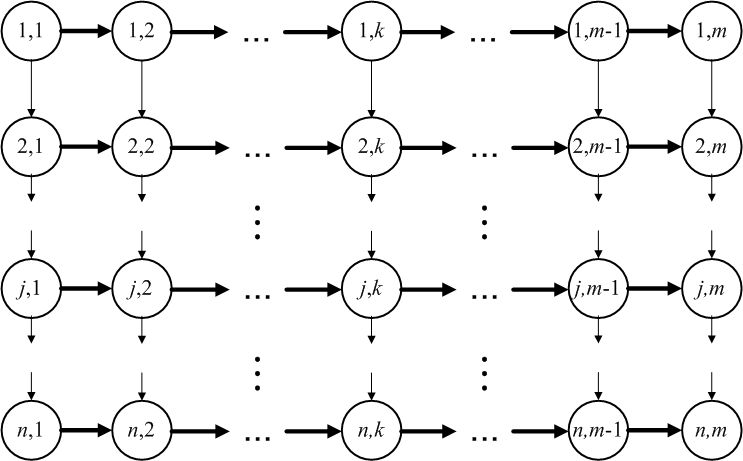}
 	\caption{Example of an expanded graph for a flow shop problem.}
 	\label{fig2}
 \end{figure}
Each vertex of such a graph can be identified with a pair (job number, machine number). 
The vertical arcs define the precedence relationships between the jobs for a particular machine. The horizontal arcs define precedence relationships between the machines for a fixed job number.
This graph is acyclic and has one start and one end vertex. 

The following assumptions hold for the FSP:
\begin{enumerate} \label{Ut}
	\item 
	Assumptions about the jobs:
	\begin{enumerate}
		\item \label{o1}
		any job cannot be performed by more than one machine at the same time;  
		\item \label{o6}  {\bf each job must be performed on each machine} according to the given machine precedences;
		\item 	any job is performed  by the machine in the shortest possible time, depending on the order of the machines;
		\item 	all jobs are equally important, meaning that there are no priorities, deadlines, or urgent jobs.
	\end{enumerate}
	\item Assumptions about the machines:
	\begin{enumerate}
		\item 
		no machine can perform more than one job at a time;
		\item   \label{o4}	after the start of a processing, each job must be completed without interruption;
		\item \label{o2}  no job is performed more than once by any machine.
	\end{enumerate}
	\item \label{p3}Assumptions about the processing times:
	\begin{enumerate}
		\item 
		the processing time of each job on each machine is independent of the sequence in which the jobs are performed;
		\item 	the processing time of each job on each machine is specified; 
		\item \label{o3}	the transport time between machines and the {\bf  setup time}, if any, {\bf are included in the job processing times};
		\item \label{o5} at least one job is being performed at any given time ({\bf no lunch breaks or overnight work stoppages}).
	\end{enumerate}
\end{enumerate}
The considered assumptions for the FSP are generally accepted and define a class of problems in the scheduling theory. 
However, in practical applications, they do not allow one to describe many types of problem statements that are widely used in practice. 
Some of them will be discussed in the next section.

The Permutation FSP (PFSP) is defined as a set of problems with such a restriction that all machines perform the jobs in the same order.
For this class of problems, an optimal schedule is constructed by changing the order in which the jobs are submitted for execution.
\section{Representing the PFSP using a recursive function}
Next, we will consider such a subset of recursive functions, the input and output of which are numerical values.
 Recursive functions have {\em arguments} and {\em parameters} as data.
The arguments of the functions are a pair $(j,k)$, where $j$ is the job number, and $k$ is the machine number.
The value of the function is the completion time of the job $j$ on the machine $k$ with a certain fixed order of jobs.
Arguments are used as arguments to a function in the usual mathematical sense.
The set of parameters ${\cal P}$ is defined outside the function and, based on the theory, can be considered as a set of externally defined recursive functions.
In fact, they are {\em global variables or constants} in terms of imperative programming.
${\cal P}=\{n,\,m,\,P(n,m)\}$ is an example of a set of parameters.     

Next, we describe the general format and properties of recursive functions that we will consider.
Let ${\cal R}$ be a recursive function with a set of arguments defined by the Cartesian product ${\cal J}\times {\cal M}$, having a set of parameters ${\cal P}$ and a set of values $N_0\cup\{\sharp\}$ (here $N_0=\{0,1,2,\dots\}$).
We will write such a function as
$${\cal R}\{{\cal P}\}:{\cal   J}\times {\cal   M}\to N_0 \cup \{\sharp\}.$$
The value $\sharp$ will indicate an infeasible schedule.
We describe the general representation of the recursive function for the PFSP in the following form:
\begin{equation}\label{f5}
%\begin{align}
%\label {op}
{\cal R}(j,k)=\left \{\begin{array}{ll}
W_1;&h_{11},h_{12},\dots, h_{1l_1};\\
W_2;&h_{21},h_{22},\dots, h_{2l_2};\\
\dots&\\
W_q;&h_{q1},h_{q2},\dots, h_{ql_q}.
\end{array}\right.
%\end{align}
\end{equation}
where $W_i$ is some expression, $h_{il}$ is a predicate, 
$q$ is the number of variants for calculating the value of 
$\cal R$.
This construction must be interpreted as a sequence:\\
$$\mbox{\it if  }h_{11}\&h_{12}\&\dots\& h_{1l_1}\mbox{\it  ~then  }{\cal R}(j,k)=W_1;$$
$$\mbox{\it if  }h_{21}\&h_{22}\&\dots\& h_{2l_2}\mbox{\it  ~then  }{\cal R}(j,k)=W_2;$$
$$\dots$$
$$\mbox{\it if  }h_{q1}\&h_{q2}\&\dots\& h_{ql_q}\mbox{\it ~then  }{\cal R}(j,k)=W_q.$$

Let $v$ and $w$ be the numbers of the right-hand sides of the function (\ref{f5}), $1\le v,w \le q$.
Then the expressions $W_1,\,W_2,\,\dots,\, W_q$ and the predicates $h_1,\,h_2,\,\dots,\,h_q$ must satisfy the following requirements:
\begin{itemize}
	\item 
The predicate sets $h_{v1}\&h_{v2}\&\dots\& h_{vl_v}$ are mutually exclusive in the sense that two different predicate sets cannot simultaneously take the value $true$, i.e.
\begin{eqnarray*}
\lefteqn{ 	\forall\,v,w\in \{1,\,2,\,\dots,\,q\},\,v\ne w\,\Rightarrow }\\
& & \Rightarrow(h_{v1}\&h_{v2}\&\dots\& h_{vl_v})\&(h_{w1}\&h_{w2}\&\dots\& h_{wl_w})=false.
\end{eqnarray*}
 \item One of all sets of predicates must necessarily take the value $true$, otherwise the function ${\cal R}(j,k)$ will not be defined everywhere, i.e.
 %\begin{equation*}
 	$$\exists v\in \{1,\,2,\,\dots,\,q\}\,\Rightarrow\,h_{v1}\&h_{v2}\&\dots\& h_{vl_v}=true.$$
 %\end{equation*}
 \item 
 In the definition of the function ${\cal R}(j,k)$, another recursive function ${\cal R}'(j,k)$ can be included in $W_v$ and $h_v$ if the chain of superpositions of the function ${\cal R}(j,k)$ does not contain a chain of calls ${\cal R}'(j,k) \Longrightarrow\dots\Longrightarrow {\cal R}(j,k)$.
\end{itemize}
{\bf Remark}.
%\begin{Remark}
{\em By default, we will assume that the recursive functions in question are defined everywhere, unless otherwise specified.}\\
%\end{Remark}
{\bf Remark}.
%\begin{remark}
{\em  Since the function ${\cal R}(j,k)$ is defined everywhere, in the case when there is no feasible schedule, the function will have the value} ${\bf\sharp}$.
%\end{remark}

In the future, instead of the term ''recursive function'', we will use the term ''function'' if there is no ambiguity.
With a fixed order of jobs, the recursive PFSP function can be defined as the mapping \cite{KLRW}
	$$C\{n,m,P(n,m)\}:{\cal J}\times {\cal M}\to  N_0\cup\{\sharp\}.$$
It looks as follows:
\begin{equation}\label{f3a}
C(j,k)=\left\lbrace \begin{array}{ll}
p_{j,k}; & j=1,\,k=1;\\
C(j,k-1)+p_{j,k}; & j=1,\,1<k\le m;\\
C(j-1,k)+p_{j,k}; & 1< j\le n,\,k=1;\\
\max\{C(j-1,k),C(j,k-1)\}+p_{j,k}; & 1< j\le n,\,1<k\le m.
\end{array}\right.
\end{equation}
Here, the value of the function $C(j,k)$ is the minimum completion time of the $j$th job on the $k$th machine with a fixed order of the jobs.
Thus, the function is defined everywhere, and the set of values $C(j,k)$ for $1\le j\le n$ and $1\le k\le m$ uniquely determines one of the optimal schedules for a fixed order of jobs. 
Each vertex $(j,k)$ of the graph in Figure \ref{fig2} corresponds to a call to the recursive function $C(j,k)$.
 \begin{figure}[!ht]
	\centering
	\includegraphics [scale=0.5]{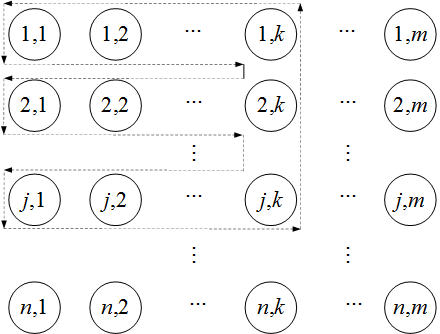}
	\caption{Illustration of vertex traversal when calculating the recursive function $C(j,k)$ in the flow shop problem.}
	\label{fig3a}
\end{figure}
The calculation of the function (\ref{f3a}) determines the vertex traversal order shown in Figure \ref{fig3a}, which we will consider.
The rather long expression $\max\{C(j-1,k),C(j,k-1)\}$ will be used frequently in the future, so for the sake of brevity, we define the operator 
\begin{equation}%\label{f3a}
\circ C(j,k)=\left\lbrace \begin{array}{ll}
%p_{j,k}; & j=1,\,k=1;\\
C(j,k-1); & j=1,\,1<k\le m;\\
C(j-1,k); & 1< j\le n,\,k=1;\\
\max\{C(j-1,k),C(j,k-1)\}; & 1< j\le n,\,1<k\le m.
\end{array}\right.
\end{equation}
%$$\circ C(j,k)=\max\{C(j-1,k),C(j,k-1)\},\; 1< j \le n,\, 1< k \le m.$$
In this case, the function (\ref{f3a}) can be written as
\begin{equation}\label{f3}
	C(j,k)=\left\lbrace \begin{array}{ll}
		p_{j,k}; & j=1,\,k=1;\\
		%C(j,k-1)+p_{j,k}; & j=1,\,1<k\le m;\\
		%C(j-1,k)+p_{j,k}; & 1< j\le n,\,k=1;\\
		\circ C(j,k)+p_{j,k}; & 1\le j\le n,\,1\le k\le m.\\
	\end{array}\right.
\end{equation}
The main optimization problem in the PFSP is to find the order of jobs in which some objective function takes the optimal value.
Most often, the objective function is the {\em makespan} (the total time it takes to complete all jobs).
In order to solve the PFSP, a function is defined (see \cite{KLRW})
$$\hat{C}\{n,m,P(n,m),{\bar\pi} \}:{\cal \hat{J}}\times {\cal M}\to N_0 \cup\{\sharp\},$$ where ${\cal \hat{J}}$ is the set of job position numbers, and the vector $\pi$ maps the set of job positions (the order in which they are started) to the set of job numbers.
\begin{equation}\label{f4}
{\hat C}(\alpha,k)=\left\lbrace \begin{array}{ll}
p_{\pi(\alpha),k}; & \alpha=1,\, k=1; \\
%{\hat C}(\alpha,k-1)+p_{\pi(\alpha),k}; & \alpha=1,\, 1<k\le m;\\
%{\hat C}(\alpha-1,k)+p_{\pi(\alpha),k}; & 1<\alpha\le n,\, k=1;\\
\circ {\hat C}(\alpha,k)+p_{\pi(\alpha),k}; &\, 1\le \alpha\le n,\, 1\le k\le m.
\end{array}\right.
\end{equation}
Here $\alpha, \alpha=1,2,...,n$ is the sequential number of the job.
Note that, when making the calculation  for all values of $1\le j\le n$ and $1\le k\le m$, the functions $C(j,k)$ and ${\hat C}(\alpha,k)$ calculate one of the optimal schedules. 

\section{Examples of recursive functions for extending flow shop problems}
The use of recursive functions makes it possible to significantly expand the PFSP statements and their importance for applications. 
Next, as an example, we will provide recursive functions for the following list of problem extensions:
\begin{enumerate}
	\item 
	time limit for receiving a job for the processing;
	\item 
periodic equipment adjustments;
	\item 
initial equipment setup;
	\item 
job sequence-dependent setup times (SDST);
	\item 
interruption of processing for a specified time period at a specified point in time;
	\item 
	time limit for completing the job ({\em deadline}).
\end{enumerate}
Problem number 7 will be considered to be the classical PFSP.
We will consider not only the extensions listed above, but also some reasonable combinations of them.
Most of the extensions under consideration are essentially additional restrictions, but for the sake of uniformity, we will use only the term ''extension''.
The proposed extensions are quite simple and therefore, they are convenient to be considered as examples. 
On the other hand, all of them have practical applications, but most of them are not reflected yet in scientific papers.
It should be taken into account that the solution of the PFSP involves rearranging jobs, which complicates the implementation of extensions.
Each extension will be a set of recursive functions and their superpositions.
The superposition of functions is a method of constructing a new function $h(x)=f(g(x))$, in which the result of one function (in this case $g$) is used as an argument for another one (function $f$).
This paper makes frequent use of a superposition.
In order to avoid describing each time all three functions $h,\,f$ and $g$ participating in the superposition, we will define that the function $h$ is the {\em superposition on the basis} of the function $g$.
In this case, it will be implied that there exists a function $f$ which, if necessary, is uniquely determined from the functions $h$ and $g$, but in this case it is not of  interest for us.\\
{\bf Remark}.
%\begin{remark}
	Newly defined recursive functions other than the initial $C(j,k)$ will be defined by a subscript, for example, $C_x(j,k)$.
%\end{remark}

The PFSP statement can be expanded by constructing a single ''universal function''.
Let there be an initial function $C(j,k)$ describing the PFSP. 
If we want to implement the first extension of the original problem statement, then this will lead to the construction of a new function $C_1(j,k)$.
If we want to implement the second extension based on the first one, this will lead to the construction of a new function $C_2(j,k)$ based on the function $C_1(j,k)$, etc.
Considering that we are going to describe six extensions just for example, the function $C(j,k)$ will become more complicated and lose clarity. 
If we assume that we will superpose all functions, then for six functions the maximum number of possible orders of superpositions is equal to the number of combinations $2^6-1=63$. 

In this paper, two methods of superposition of recursive functions are used. 
The first method is used when calculations can and should be associated with a machine.
In this case, each machine is assigned a type attribute, and each type has its own recursive function associated with it.
This is possible when the function calculates the completion time of the $j$th job on the $k$th machine.
However, this method is inconvenient in the case when the execution of a recursive function must be associated with each machine, for example, when checking the deadline constraint.
The second method is to superpose some function on the basis of each function associated with the machine.

 \begin{figure}[!ht]
	\centering
	\includegraphics [scale=0.6]{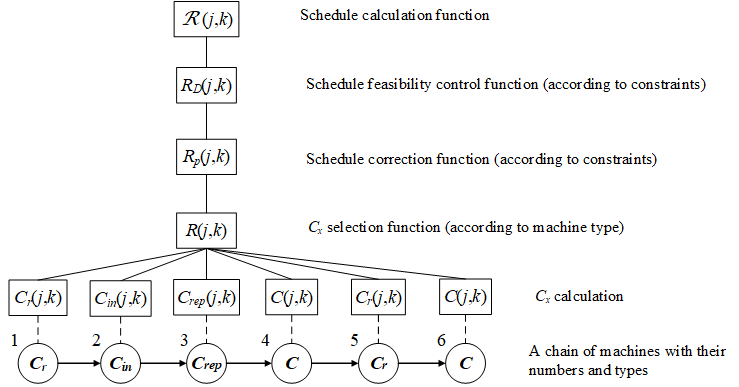}
	\caption{An example of the structure of superpositions of PFSP functions}
	\label{fig17}
\end{figure}
Let us consider the general scheme of a superposition using the example in Figure \ref{fig17}. 
At the upper level, the function ${\cal R}(j,k)$ is used to calculate the final completion time of the $j$th job on the $k$th machine.
At the second level, the function $R_D(j,k)$ controls the allowable completion time in accordance with the time constraint.
At the third level, using the function $R_p(j,k)$, the completion time is adjusted according to the time constraint. 
At the fourth level, using the function $R(j,k)$, a specific recursive function is selected according to the type of the machine.
At the fifth level, using a function associated with a machine of some type, the completion time of the $j$th job on the $k$th machine is calculated.
From the structure considered, it can be seen that the functions of schedule correction and schedule feasibilty control are common to all machines.
The remaining functions correspond to a specific machine. Next, let us look at the specific application of these methods.

At the first stage, we will describe a set of functions associated with machines.
To describe the specifics of the extended tasks, we define for the machines in Figure  \ref{fig1a} the set of {\em types of machines} $\Im=\{\bf{C},\,\bf{C_r},\,\bf{C_{rep}},\,\bf{C_{in}},\,\bf{C_{SDST}} \}$. 
The types of different machines can match, i.e., the number of types of machines is less than or equal to $m$.
For each type of machine ${\bf C_x}$, we will define a recursive function $C_x(j,k)$. 
In accordance with this principle, we define the set of recursive functions $$\{C(j,k),\,C_r(j,k),\,C_{rep}(j,k),\,C_{in}(j,k),\,C_{SDST}(j,k)\}$$ and the $type:\cal M\to  \Im$, which uniquely assigns a certain type to each machine.
The function $R(j,k)$ activates the corresponding recursive function as follows:\\
For a fixed order of jobs:
\begin{equation}\label{r1}
	R(j,k)=\left\lbrace \begin{array}{ll}
C(j,k); & type(k)={\bf C};\\
C_r(j,k); &  type(k)={\bf C_r};\\
%\qquad\dots&\\
%C_D(j,k); &  type(k)={\bf C_D};\\
C_{rep}(j,k); & type(k)={\bf C_{rep}};\\
C_{in}(j,k); & type(k)={\bf C_{in}};\\
C_{SDST}(j,k); & type(k)={\bf C_{SDST}}.\\
%C_{p}(j,k); & type(k)={\bf C_p}.
\end{array}\right.\\
\end{equation}
For the PFSP:
\begin{equation}\label{r2}
{\hat R}(j,k)=\left\lbrace \begin{array}{ll}
{\hat C}(j,k); & type(k)={\bf C};\\
{\hat C}_r(j,k); &  type(k)={\bf C_r};\\
%\qquad\dots&\\
%{\hat C}_D(j,k); &  type(k)={\bf C_D};\\
{\hat C}_{rep}(j,k); & type(k)={\bf C_{rep}};\\
{\hat C}_{in}(j,k); & type(k)={\bf C_{in}};\\
{\hat C}_{SDST}(j,k); & type(k)={\bf C_{SDST}}.\\
%{\hat C}_{p}(j,k); & type(k)={\bf C_p}.
\end{array}\right.
\end{equation}
In this case, a call of the function $R$ or $\hat R$ is transformed into a call of some specific function, depending on the machine type.  
First, we will consider recursive functions for a fixed order of jobs.
This makes the function easier to describe and understand.
Then examples will be given for these functions.
The permutation version of the recursive function will be described after all the explanations.
An important requirement for these functions is the correctness of their calculation when the order of execution of the jobs changes.
We will not write about this further for each function, assuming this is a default requirement. 

{\bf PFSP}. The first function $C(j,k)$ from the right side of the formula (\ref{r1}) is associated with a machine of type ${\bf C}$ and implements the calculation of the job completion time for the PFSP:
$$C\{n,m,P(n,m)\}:{\cal J}\times {\cal M}\to  N_0\cup\{\sharp\}.$$
\begin{equation}\label{f3b}
C(j,k)=\left\lbrace \begin{array}{ll}
p_{j,k}; & j=1,\,k=1;\\
%{\cal R}(j,k-1)+p_{j,k}; & j=1,\,1<k\le m;\\
%{\cal R}(j-1,k)+p_{j,k}; & 1< j\le n,\,k=1;\\
\circ {\cal R}(j,k)+p_{j,k}; & 1\le j\le n,\,1\le k\le m.
\end{array}\right.
\end{equation}
Formula (\ref{f3b}) is based on formula (\ref{f3}), in which, on the right side, the function $C(j,k)$ is replaced by the function ${\cal R}(j,k)$ (see the function ${\cal R}(j,k)$ in Figure  \ref{fig17}). 
In the example in Figure  \ref{fig17}, there are three machines of type $C$: 2, 4 and 6.\\
Permutation version:\\
 $$\hat{C}\{n,m,P(n,m),{\bar\pi} \}:{\cal \hat{J}}\times {\cal M}\to N_0 \cup\{\sharp\},$$
 \begin{equation} \label{f4p}
 {\hat C}(\alpha,k)=\left\lbrace \begin{array}{ll}
 p_{\pi(\alpha),k}; & \alpha=1,\, k=1; \\
%{\hat {\cal R}}(\alpha,k-1)+p_{\pi(\alpha),k}; & \alpha=1,\, 1<k\le m;\\
%{\hat {\cal R}}(\alpha-1,k)+p_{\pi(\alpha),k}; & 1<\alpha\le n,\, k=1;\\
 \circ {\hat {\cal R}}(\alpha,k)+p_{\pi(\alpha),k}; &\, 1\le\alpha\le n,\, 1\le k\le m.
 \end{array}\right.
 \end{equation}
 
 {\bf Constraint $r_j$}.
The next function from the right side of formula (\ref{r1}) $C_r(j,k)$ is associated with a machine of type ${\bf C_r}$, which must be the first in the chain of machines.
 This function implements a limitation on the start time of the $j$-th job. 
The  parameter $\bar{r}$ determines the point in time from which the job can be put into service, but it is not guaranteed that its processing will begin at that point.
The constraint can be defined for each job by the vector $\bar{r}\;(|\bar{r}|=n)$. 
Since each job starts on the first machine, the function can be represented as follows:
	$$C_r\{n,m,\bar{r},P(n,m)\}:{\cal J}\times {\cal M}\to N_0\cup \{\sharp\}.$$
\begin{equation}\label{f3c}
	C_r(j,k)=\left\lbrace \begin{array}{ll}
		r_j+p_{j,k}; &   j=1,\,k=1;\\
		%{\cal R}(j,k-1)+p_{j,k}; & j=1,\,1<k\le m;\\
		\max \{r_j,\circ{\cal R}(j,k)\}+p_{j,k}; & j=1,\,1<k\le m;\\
		\circ {\cal R}(j,k)+p_{j,k}; & 1< j\le n,\,1\le k\le m.
	\end{array}\right.
\end{equation}
From the description of the function $C_r(j,k)$, it is clear that the changes in formula (\ref{f3b}) concern the first machine, but with the help of function (\ref{f3c}), the correct values are calculated for a machine of type ${\bf C}$ or ${\bf C_r}$ in any position. 

\begin{figure}[!ht]
	\centering
	\includegraphics [scale=0.6]{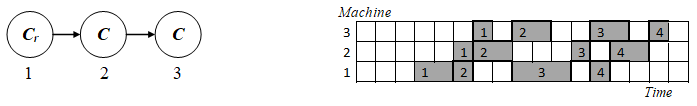}
	\caption{ An example of implementing the $r_j$ constraint.}
\label{fig15}
\end{figure}

{\bf Example 1.} 
%\begin{Example}\label{ex1}
	Let us consider an example with $n=4,\,m=3$ and the vector $r=(3,3,8,12)$. The matrix of the processing times $p_{j,k}$ is as follows:
$$
%		\begin{equation*}
		P=\left( \begin{array}{ccc}
			2&1&1 \\
			1&2&2\\ 
			3&1&2 \\
			1&2&1\end{array} \right).
%	\end{equation*}
$$
The first machine in the chain will be a machine of type ${\bf C_r}$. 
The remaining machines are of type ${\bf C}$, i.e., they correspond to the recursive function (\ref{f3b}).
The chain of machines and the timing diagram are shown in Figure  \ref{fig15}. 
%\end{Example}

Permutation version:\\
%----------------------------------------------------
$${\hat C}_r\{n,m,\bar{r},P(n,m),\bar \pi\}:{\cal J}\times {\cal M}\to N_0\cup \{\sharp\}.$$
$$
%\begin{equation*}\label{f3p}
	{\hat C}_r(\alpha,k)=\left\lbrace \begin{array}{ll}
		r_{\pi(\alpha)}+p_{\pi(\alpha),k}; &   \alpha=1,\,k=1;\\
		%{\hat {\cal R}}(\alpha,k-1)+p_{\pi(\alpha),k}; & \alpha=1,\,1<k\le m;\\
		\max \{r_{\pi(\alpha)},\circ{\hat {\cal R}}(\alpha,k)\}+p_{\pi(\alpha),k}; &  \alpha=1,\,1<k\le m;\\
		\circ {\hat {\cal R}}(\alpha,k)+p_{\pi(\alpha),k}; & 1< \alpha\le n,\,1\le k\le m.
	\end{array}\right.
%\end{equation*}
$$
{\bf Performing periodic equipment adjustments}.
A common manufacturing operation is the periodic change of tools, for example, to sharpen drills or change the equipment after completing $q$ jobs.
This means that a machine of type ${\bf C_{rep}}$ with number $k$ processes each job with number $j$ not a multiple of $q$, i.e. $j\,mod\,q\ne0$, in time $p_{j,k}$, and jobs with a number multiple of $q$ in time $p_{j,k}+\stackrel{rep}{\Delta}_k$.
Here $j\,mod\,q$ is the remainder of dividing integers $j$ by $q$. 
$\stackrel{rep}{\Delta}$ has the index $k$ because there can be several such machines with adjustments in the chain, and each of them will have its own value of $\stackrel{rep}{\Delta}$.
$$C_{rep}\{n,m,q,\|\stackrel{rep}{\Delta}_k\|,P(n,m)\}:{\cal J}\times {\cal M}\to  N_0\cup\{\sharp\}.$$
Here $q>1$ is the number of jobs after which the machine performs equipment adjustments; $\stackrel{rep}{\Delta}_k$ is the time required to perform periodic adjustments
(equipment retooling) of machine $k$:
\begin{equation} 	\label {rep}
C_{rep}(j,k)=\left \{\begin{array}{ll}
p_{j,k}; &j=1,\,k=1;\\
%{\cal R}(j-1,k)+p_{j,k}; &1< j \le n,\,k=1,\,j\,mod\,q\ne 0;\\
%{\cal R}(j-1,k)+p_{j,k}+\stackrel{rep}{\Delta}_k; &1< j \le n,\,k=1,\,j\,mod\,q=0;\\
%{\cal R}(j,k-1)+p_{j,k}; &j=1,\,1< k\le m;\\	
\circ {\cal R}(j,k)+p_{j,k}; &1\le  j \le n,\,1\le k\le m,\,j\,mod\,q\ne 0;\\
\circ {\cal R}(j,k)+p_{j,k}+\stackrel{rep}{\Delta}_k; &1\le j \le n,\,1\le k\le m,\,j\,mod\,q=0.\\
\end{array}\right.
\end{equation}
From the definition of the function, it is clear that it is correctly executed for any number of machines of type ${\bf C_{rep}}$.
Each machine of type ${\bf C_{rep}}$ can have its own $\stackrel{rep}{\Delta}_k$ value.

\begin{figure}[!ht]
	\centering
	\includegraphics [scale=0.8]{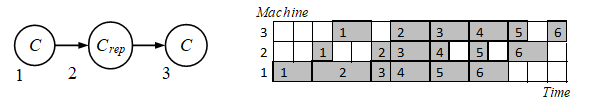}
	\caption{An example of a timing diagram for a periodic equipment adjustment on machine 2 with period 3.}
	\label{fig14}
\end{figure}

{\bf Example 2.}  
%\begin{Example}
Figure \ref{fig14} shows an example of a graph and the timing diagram for three machines and six jobs.
The second machine performs the main work and periodically adjusts the equipment. 
We have $q=3,\, \stackrel{rep}{\Delta}=1$ (without index, since there is only one machine with the adjustment).
The matrix $||p_{j,k}||$ is as follows:
$$
%	\begin{equation*}
	P=\left( \begin{array}{ccc}
		2&1&2\\
		3&1&2\\ 
		1&1&2 \\
		2&1&2\\
		2&1&1\\
		2&1&1\\
	\end{array} \right).
%\end{equation*}
$$
From the diagram, it is clear that machine 2 performs the jobs 1, 2, 4, 5 in time 1, and jobs 3 and 6 in time 2.
It might seem that it is sufficient to simply adjust the corresponding times for machine 2 in the matrix $||p_{j,k}||$.
However, in this case, it will be necessary to adjust the matrix each time the order of work is changed.
%\end{Example}

Permutation version:
%---------------------------------------------------------
$${\hat C}_{rep}\{n,m,q,\|\stackrel{rep}{\Delta}_k\|,P(n,m),{\bar \pi}\}:{\cal J}\times {\cal M}\to  N_0\cup\{\sharp\}.$$
$$
%\begin{equation*} 	\label {repp}
	{\hat C}_{rep}(\alpha,k)=\left \{\begin{array}{ll}
		p_{\pi(\alpha),k}; &\alpha=1,\,k=1;\\
		%{\hat {\cal R}}(\alpha-1,k)+p_{\pi(\alpha),k}; &1< \alpha \le n,\,k=1,\,
        %\alpha\,mod\,q\ne 0;\\
		%{\hat {\cal R}}(\alpha-1,k)+p_{\pi(\alpha),k}+\stackrel{rep}{\Delta}_k; &1< \alpha %\le n,\,k=1,\,\alpha\,mod\,q=0;\\
		%{\hat {\cal R}}(\alpha,k-1)+p_{\pi(\alpha),k}; &\alpha=1,\,1< k\le m;\\	
		\circ {\hat {\cal R}}(\alpha,k)+p_{\pi(\alpha),k}; &1\le \alpha \le n,\,1\le k\le m,\,
        \alpha\,mod\,q\ne 0;\\
		\circ {\hat {\cal R}}(\alpha,k)+p_{\pi(\alpha),k}+\stackrel{rep}{\Delta}_k; &1\le \alpha \le n,\,1\le k\le m,\,
        \alpha\,mod\,q=0.\\
	\end{array}\right.
%\end{equation*}
$$
{\bf Initial equipment setup}.
Sometimes the first job on a particular machine takes longer to complete than subsequent jobs.
This is related to the initial setup of the equipment.
This type of setup work is performed once, and then its execution time can be considered equal to zero. 
\begin{figure}[!ht]
	\centering
	\includegraphics [scale=0.56]{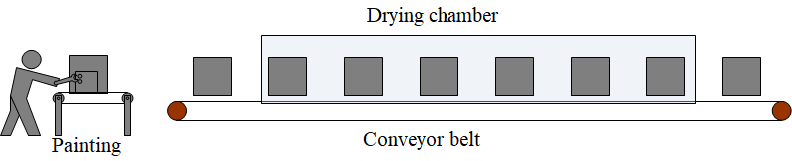}
	\caption{An example of performing work with the initial setup.}
	\label{fig11}
\end{figure}
An additional time may be due to the technological process.
Figure \ref{fig11} shows a drying chamber that dries parts after painting.
Each part lying on the conveyor belt passes through the chamber for a certain time $T$, during which the part dries.
In this case, the first part will exit the chamber after time $T$, and the subsequent ones will exit after time $T/6$, since there are 6 parts in the chamber at the same time. 
The setup can be implemented as a machine of type ${\bf C_{in}}$, which, when servicing the first job, performs the setup job in time $\stackrel{in}{\Delta}_{j,k}$, and then the main job in time $p_{1,k}$.
$\stackrel{in}{\Delta}_{j,k}$ is an element of the matrix $\stackrel{in}{\Delta}$.
The index $k$ is used because such a machine may not be the first in the chain, and there may be several machines performing the initial setup.
The index $j$ determines the job number for the case, where the $j$-th job becomes the first one in the permutation of the jobs. 
The function corresponding to the machine with the initial setup can be described as follows:\\
$$C_{in}\{n,m,P(n,m),\|\stackrel{in}{\Delta}_{j,k}\| \}:{\cal J}\times {\cal M}\to N_0\cup\{\sharp\}, $$
\begin{equation} 	\label {in}
C_{in}(j,k)=\left \{\begin{array}{ll}
\stackrel{in}{\Delta}_{j,k}+p_{j,k}; &j=1,\,k=1;\\
\stackrel{in}{\Delta}_{j,k}+{\cal R}(j,k-1)+p_{j,k}; &j=1,\,1< k\le m;\\
%{\cal R}(j-1,k)+p_{j,k}; &	1<j\le n,\,k=1;\\
\circ {\cal R}(j,k)+p_{j,k}; &1< j\le n,\,1\le k\le m;\\
\end{array}\right.
\end{equation}

\begin{figure}[!ht]
	\centering
	\includegraphics [scale=0.68]{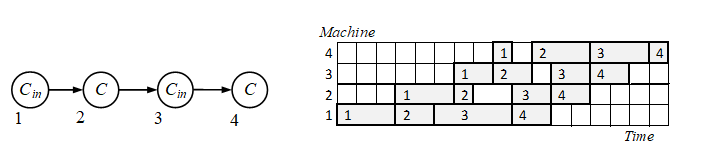}
	\caption{Examples of using the initial equipment setup function.}
	\label{fig10}
\end{figure}

{\bf Example 3.}  
%\begin{Example}
	Figure \ref{fig10} shows a chain of machines and a diagram for the order of jobs 1, 2, 3, 4.
  The matrices $P$ and $\stackrel{in}{\Delta}$  are as  follows:
  $$
%	\begin{equation*}
		P=\left( \begin{array}{cccc}
			1&3&1&1\\
			2&1&2&3\\ 
			4&2&2&3 \\
			2&2&2&1\\
		\end{array} \right),\qquad
	\stackrel{in}{\Delta}=\left( \begin{array}{cccc}
	2&0&1&0 \\
	0&0&0&0\\ 
	0&0&0&0 \\
	0&0&0&0\end{array} \right).
%	\end{equation*}
$$
The matrix $\stackrel{in}{\Delta}$ has non-zero values only in the first row because an example without job permutations is considered.
%\end{Example}

Permutation version:
%---------------------------------------------------------
$${\hat C}_{in}\{n,m,P(n,m),\|\stackrel{in}{\Delta}_{j,k}\|,{\bar \pi} \}:{\cal J}\times {\cal M}\to N_0\cup\{\sharp\}, $$
$$
% \begin{equation*} 	\label {inp}
{\hat C}_{in}(\alpha,k)=\left \{\begin{array}{ll}
\stackrel{in}{\Delta}_{\pi(\alpha),k}+p_{\pi(\alpha),k}; &\alpha=1,\,k=1;\\
\stackrel{in}{\Delta}_{\pi(\alpha),k}+{\hat {\cal R}}(\alpha,k-1)+p_{\pi(\alpha),k}; &\alpha=1,\,1< k\le m;\\
%{\hat {\cal R}}(\alpha-1,k)+p_{\pi(\alpha),k}; &	1<\alpha\le n,\,k=1;\\
\circ {\hat {\cal R}}(\alpha,k)+p_{\pi(\alpha),k}; &1<\alpha\le n,\,1\le k\le m.\\
\end{array}\right.
%\end{equation*}
$$
{\bf Job sequence-dependent setup times (SDST).}
The problem of planning an SDST-type production process is one of the frequently encountered problems in the scheduling theory \cite{Sad,Sat,Ci}.
Here, the setup time depends on the sequence of jobs.
This situation usually occurs when an enterprise uses multi-purpose equipment that needs to be configured.
Two examples of SDST workflow planning are given in \cite{Sat}:
\begin{enumerate}
	\item 
The textile industry, where setup times are significant, as fabric types are assigned to weaving machines equipped with winding chains.
When changing the type of fabric on the machine, it is necessary to replace the winding chain.
The time it takes to do this depends on the previous and current fabric types.
\item 
Stamping machines used by most car manufacturers.
On such machines, the setup time between  production of parts depends on the sequence of work, as it is associated with the replacement of heavy dies.
\end{enumerate}
In \cite{Sad}, the setup times of the $k$-th $(k=1,2,\dots, m)$ machine are defined for $n^2$ pairs of {\em previous job}$\times${\em current job}.
In reality, the number of combinations is determined not by the number of jobs, but by the number of tools, which are significantly less.
We will assume that the function $Q:{\cal J}\to Z$ maps the set of jobs to the set of tools, and $z=|Z|$ is the number of tools.
The setup times of the $k$-th machine for a pair of $i,j$ jobs are determined by the array $\tau_{v,w,k}$ of dimension $z\times z\times m$ and $v=Q(i),\, w=Q(j)$.
The values of $\tau_{Q(j),Q(j),k},\; 1\le j\le n$, will determine the initial setup times of the equipment when the job is the first in the schedule.
The function will be described as follows:
$$C_{SDST}\{n,m,z,P(n,m),\|\tau_{v,w,k}\|\}:{\cal J}\times {\cal M}\to N_0\cup  \{\sharp\}.$$
\begin{equation}\label{sdstj}
C_{SDST}(j,k)=\left\lbrace \begin{array}{ll}
p_{j,k}+\tau_{Q(j),Q(j),k}; & j=1,\, k=1; \\
%{\cal R}(j,k-1)+p_{j,k}+\tau_{Q(j),Q(j),k}; & j=1,\, 1<k\le m;\\
%{\cal R}(j-1,k)+p_{j,k}+\tau_{Q(j-1),Q(j),k}; & 1<j \le n,\, k=1;\\
\circ {\cal R}(j,k)+p_{j,k}+\tau_{Q(j-1),Q(j),k}; &\, 1\le j \le n,\, 1\le k\le m.
\end{array}\right.
\end{equation}
Let us also describe the permutation function: 
$${\hat C}_{SDST}\{n,m,z,P(n,m),{\bar \pi},\|\tau_{v,w,k}\|\}:{\cal J}\times {\cal M}\to N_0\cup  \{\sharp\}.$$
$$
%\begin{equation*}\label{sdst}
{\hat C}_{SDST}(\alpha,k)=
        \left\lbrace \begin{array}{ll}
p_{\pi(\alpha),k}+\tau_{Q(\pi(\alpha)),Q(\pi(\alpha)),k}; & \alpha=1,\, k=1; \\
%{\cal {\hat R}}(\alpha,k-1)+p_{\pi(\alpha),k}+\tau_{Q(\pi(\alpha)),Q(\pi(\alpha)),k}; & \alpha=1,\, 1<k\le m;\\
%{\cal {\hat R}}(\alpha-1,k)+p_{\pi(\alpha),k}+\tau_{Q(\pi(\alpha-1)),Q(\pi(\alpha)),k}; & 1<\alpha\le n,\, k=1;\\
\circ {\cal {\hat R}}(\alpha,k)+p_{\pi(\alpha),k}+\tau_{Q(\pi(\alpha-1)),Q(\pi(\alpha)),k}; &\, 1\le \alpha\le n,\, 1\le k\le m.
\end{array}\right.
%\end{equation*}
$$
Here $\pi=(\alpha_1,\,\alpha_2,\dots,\alpha_n),\; Q(\pi(\alpha))$ is the number of the tool for performing the job with the ordinal number $\alpha,\; \alpha=1,\,2,\dots,n$, and $Q(\pi(\alpha-1))$ is the number of the tool for performing the previous job with the ordinal number $\alpha-1$.

Next, in this section, we will consider second- and third-level recursive functions (see Figure \ref{fig17}). These functions must implement the correction and execution time control for all $(j,k)$ pairs. 

{\bf Processing interruption}.
Real production processes are necessarily interrupted for a certain fixed period of time.
For example, there might be a lunch break, a break between shifts, a break at night before the next shift, etc.
If the execution of jobs is interrupted (violation of the assumption on page \pageref{p3}, property \ref{o5}), then the function implementing this feature looks as follows:
$$R_p\{n,m,T_{e},\stackrel{p}{\Delta},P(n,m)\}:{\cal J}\times {\cal M}\to  N_0\cup \{\sharp\},$$
where $T_{e}$ is the pause start time (end of work);
$\stackrel{p}{\Delta}$ is the value of the interval for which the process stops, for example, an hour for lunch.
In this case, the function $R(j,k)$ is a result of a superposition based on the function $R_p(j,k)$ (see Figure  \ref{fig17}).

\begin{figure}[!ht]
	\centering
	\includegraphics[width=0.8\linewidth]{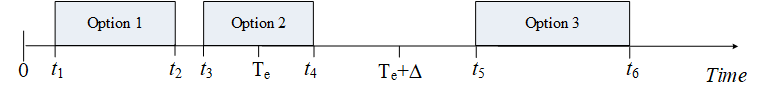}
	\caption{Job execution options.}
	\label{fig11a}
\end{figure}

This is because the pause-induced job completion time adjustment applies to all machines.
Since the description of the function is quite large, we will describe its structure.
Figure \ref{fig11a} shows possible options for executing the $j$-th job on the $k$-th machine:
\begin{enumerate}
	\item 
	The execution of the $j$-th job on the $k$-th machine starts at time $t_1$ and ends at time $t_2$ such that $t_1,\,t_2\le T_e$. The completion time is computed by function (\ref{f3}).
	\item \label{co2}
	The execution of the $j$-th job on the $k$-th machine starts at time $t_3\le T_e$ and ends at time $t_4>T_e$. In this case, the start time of the job increases to $T_e+\stackrel{p}{\Delta}$.
	\item 
	The execution of the $j$-th job on the $k$-th machine begins at time $t_5>T_e$. In this case, the completion time is calculated by function (\ref{f3}), since the start time is adjusted when the job is executed by the previous machine, i.e., there is a machine preceding the one under consideration for which the option \ref{co2} is satisfied.
\end{enumerate}
Thus, the changes to the completion time calculation algorithm only apply to option 2.
\begin{equation} \label{f4a}
	R_p(j,k)=\left\lbrace \begin{array}{lll}
		R(j,k); & R(j,k)\le T_e, &1\le j\le n,\, 1\le k\le m;\\
		%R(j,k); & {\cal R}(j,k-1)> T_e, & j=1,\, 1<k\le m;\\
        %R(j,k); & {\cal R}(j-1,k)> T_e, & 1<j\le n,\, k=1;\\
		R(j,k); & \circ {\cal R}(j,k)> T_e, &1\le j\le n,\, 1\le k\le m;\\
		% v1 -- j=1,  1<k
		%T_{e}+\stackrel{p}{\Delta}+ &{\cal R}(j,k-1)\le T_e, &\\
	%\;\;+(R(j,k)-{\cal R}(j,k-1)); &
    %R(j,k)>T_e, &j=1,\, 1< k\le m;\\
    % v2   --  1<j,  k=1
		%T_{e}+\stackrel{p}{\Delta}+ &{\cal R}(j-1,k)\le T_e, &\\
	%\;\;+(R(j,k)-{\cal R}(j-1,k)); &
    %R(j,k)>T_e, &1<j\le n,\,k=1;\\
    %  v3   -- 1<j, 1<k
		T_{e}+\stackrel{p}{\Delta}+ &\circ {\cal R}(j,k)\le T_e, &\\
	\;\;+(R(j,k)-\circ {\cal R}(j,k)); &
    R(j,k)>T_e, &1\le j\le n,\, 1\le k\le m.\\    
	\end{array}\right.
\end{equation}
This formula means that for any type of machine, the value of the function $R(j,k)$ is calculated and adjusted in accordance with the values of $T_e$ and $\stackrel{p}{\Delta}$.
The value $\circ {\cal R}(j,k)$ is used here as the start time of the $j$-th job on the $k$-th machine, since using the value ${\cal R}(j,k)-p_{j,k}$ in this case is incorrect.
The time it takes to complete a job may vary, such as in the case of a periodic job.

\begin{figure}[!ht]
		\centering
		\includegraphics [scale=0.6]{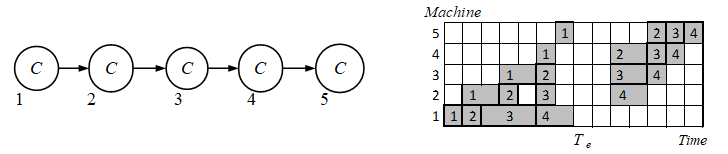}
		\caption{Example of a timing diagram for a problem with a pause.}
		\label{fig12}
	\end{figure}

{\bf Example 4}.
%\begin{Example}
	Consider an example of processing interruption with the following parameters:\\
	$n=4,\; m=5,\; T_{e}=7,\; \stackrel{p}{\Delta}=2.$\\
	The matrix $P(n,m):$ is as follows: 
$$    
%	\begin{equation*}
		P=\left( \begin{array}{ccccc}
			1&2&2&1&1 \\
			1&1&1&2&1\\ 
			3&1&2&1&1 \\
			2&2&1&1&1\end{array} \right).
%	\end{equation*}
$$
	An example diagram of the corresponding schedule for the order of tasks $(1, 2, 3, 4)$ is shown in Figure \ref{fig12}.
%\end{Example}

Permutation version:\\
%---------------------------------------------------------
$${\hat R}_p\{n,m,T_{e},\stackrel{p}{\Delta},P(n,m),{\bar \pi}\}:{\cal J}\times {\cal M}\to  N_0\cup \{\sharp\},$$
$$
%\begin{equation*} \label{f4ap}
{\hat R}_p(\alpha,k)=\left\lbrace \begin{array}{lll}
{\hat R}(\alpha,k); & {\hat R}(\alpha,k)\le T_e,&1\le \alpha\le n,\, 1\le k\le m;\\
%{\hat R}(\alpha,k); & {\hat {\cal R}}(\alpha,k-1)> T_e, & \alpha=1,\, 1<k\le m;\\
%{\hat R}(\alpha,k); & {\hat {\cal R}}(\alpha-1,k)> T_e, & 1<\alpha\le n,\, k=1;\\
{\hat R}(\alpha,k); & \circ {\hat{\cal  R}}(\alpha,k)> T_e,&1\le \alpha\le n, 1\le k\le m;\\
		% v1 -- j=1,  1<k
		%T_{e}+\stackrel{p}{\Delta}+ &{\hat {\cal  R}}(\alpha,k-1)\le T_e, &\\
	%\;\;+({\hat R}(\alpha,k)-{\hat {\cal  R}}(\alpha,k-1)); &
    %{\hat R}(\alpha,k)>T_e, &j=1,\, 1< k\le m;\\
    % v2   --  1<j,  k=1
		%T_{e}+\stackrel{p}{\Delta}+ &{\cal R}(j-1,k)\le T_e, &\\
	%\;\;+({\hat R}(\alpha,k)-{\cal R}(j-1,k)); &
   %{\hat R}(\alpha,k)>T_e, &1<j\le n,\,k=1;\\
    %  v3   -- 1<j, 1<k
T_{e}+\stackrel{p}{\Delta}+&\circ {\hat{\cal  R}}(\alpha,k)\le T_e,\\
\;\;+({\hat R}(\alpha,k)-\circ {\hat{\cal  R}}(\alpha,k)); &\,{\hat R}(\alpha,k)>T_e,&1\le \alpha\le n,\, 1\le k\le m.\\
%&\qquad\qquad\qquad\qquad  \qquad\qquad    1\le \alpha\le n,\, 1\le k\le m.\\
\end{array}\right.
%\end{equation*}
$$
{\bf Deadline limitation}. Let $D_j$ be the deadline for completing the $j$-th job. 
A schedule that contains a job that is completed after its deadline is infeasible.
Since the constraint $D_j$ must be checked for each machine, it makes sense to use the function $R_D(j,k)$ (see the function $R_D(j,k)$ in Figure \ref{fig17}).
$$R_D\{n,m,\bar{D}\}:{\cal J}\times {\cal M}\to N_0\cup \{\sharp\}.$$
This function  implements only the check of the $D_j$ constraint:
$$
%\begin{align} 	\label {Dj}
	R_{D}(j,k)=\left \{\begin{array}{ll}
		R_p(j,k); &  R_p(j,k)\le D_j,\,1\le j\le n,\,1\le  k\le m;\\			
		\sharp; & R_p(j,k)> D_j,\,1\le  j\le n,\,1\le  k\le m.\\
	\end{array}\right.
%\end{align}
$$
{\bf Example 5}. 
%\begin{Example}
Let us consider Example 1 %\ref{ex1}
and supplement it with a check of the constraint $D$.
Let $D=(8,13,15,16)$.
For the order of jobs be $(1,\,2,\,3,\,4)$, thus constraint $D$ is satisfied.
For the vector $D=(7,13,14,15)$, the constraint is not satisfied for the 4th job, i.e., this schedule is not feasible.
%\end{Example}

Permutation version:\\
%---------------------------------------------------------
$$R_D\{n,m,\bar{D},{\bar \pi}\}:{\cal J}\times {\cal M}\to N_0\cup \{\sharp\}.$$
$$
%\begin{align*} 	%\label {Dj}
	{\hat R}_{D}(\alpha,k)=\left \{\begin{array}{ll}
		{\hat R}_p(\alpha,k); &  {\hat R}_p(\alpha,k)\le D_{\pi(\alpha)},\,1\le \alpha\le n,\,1\le  k\le m;\\			
		\sharp; & {\hat R}_p(\alpha,k)> D_{\pi(\alpha)},\,1\le  \alpha\le n,\,1\le  k\le m.\\
	\end{array}\right.
%\end{align*}
$$

This chain of function superpositions will not lead to a loop of calls, since it always ends with a call of $C_x(j,k)$ to ${\cal R}(j-1,k)$
or ${\cal R}(j,k-1)$ (see formulas (\ref{f3b}), (\ref{f3c}),(\ref{rep}),(\ref{in}),(\ref{sdstj})). 

For a fixed order of jobs, the schedule is calculated as a set of values of the function 
$$
%\begin{equation*}
	{\cal R}(j,k)=R_D(j,k),\,1\le j\le n,\,1\le  k\le m.\\
%\end{equation*}
$$
for the permutation
$$
%\begin{equation*}
	{\hat {\cal R}}(\alpha,k)={\hat R}_D(\alpha,k),\,1\le\alpha\le n,\,1\le  k\le m.\\
%\end{equation*}
$$

For a fixed order of jobs, the general set of recursive functions implementing the stated PFSP extensions is as follows:\\
\begin{equation}\label{ff1}
{\cal R}(j,k)=R_D(j,k).\\
\end{equation}
$$
%\begin{align*} 	\label {ff2}
	R_{D}(j,k)=\left \{\begin{array}{ll}
		R_p(j,k); &  R_p(j,k)\le D_j,\,1\le j\le n,\,1\le  k\le m;\\			
		\sharp; & R_p(j,k)> D_j,\,1\le  j\le n,\,1\le  k\le m.\\
	\end{array}\right.
%\end{align*}
$$
$$
%\begin{equation} \label{f4a}
	R_p(j,k)=\left\lbrace \begin{array}{lll}
		R(j,k); & R(j,k)\le T_e, &1\le j\le n,\, 1\le k\le m;\\
		%R(j,k); & {\cal R}(j,k-1)> T_e, & j=1,\, 1<k\le m;\\
        %R(j,k); & {\cal R}(j-1,k)> T_e, & 1<j\le n,\, k=1;\\
		R(j,k); & \circ {\cal R}(j,k)> T_e, &1\le j\le n,\, 1\le k\le m;\\
		% v1 -- j=1,  1<k
		%T_{e}+\stackrel{p}{\Delta}+ &{\cal R}(j,k-1)\le T_e, &\\
	%\;\;+(R(j,k)-{\cal R}(j,k-1)); &
    %R(j,k)>T_e, &j=1,\, 1< k\le m;\\
    % v2   --  1<j,  k=1
		%T_{e}+\stackrel{p}{\Delta}+ &{\cal R}(j-1,k)\le T_e, &\\
	%\;\;+(R(j,k)-{\cal R}(j-1,k)); &
    %R(j,k)>T_e, &1<j\le n,\,k=1;\\
    %  v3   -- 1<j, 1<k
		T_{e}+\stackrel{p}{\Delta}+ &\circ {\cal R}(j,k)\le T_e, &\\
	\;\;+(R(j,k)-\circ {\cal R}(j,k)); &
    R(j,k)>T_e, &1\le j\le n,\, 1\le k\le m.\\    
	\end{array}\right.
%\end{equation}
$$
\begin{equation} \label{ff4}
	R(j,k)=\left\lbrace \begin{array}{ll}
		C(j,k); & type(k)={\bf C};\\		
		C_r(j,k); & type(k)={\bf Cr};\\
		C_{rep}(j,k); &  type(k)={\bf Crep};\\
		C_{SDST}(j,k); & type(k)={\bf C{\scriptstyle SDST}};\\		
		C_{in}(j,k); &  type(k)={\bf Cin}.
	\end{array}\right.
\end{equation}
\begin{equation}\label{ff5}
C(j,k)=\left\lbrace \begin{array}{ll}
p_{j,k}; & j=1,\,k=1;\\
%{\cal R}(j,k-1)+p_{j,k}; & j=1,\,1<k\le m;\\
%{\cal R}(j-1,k)+p_{j,k}; & 1< j\le n,\,k=1;\\
\circ {\cal R}(j,k)+p_{j,k}; & 1\le j\le n,\,1\le k\le m.
\end{array}\right.
\end{equation}
\begin{equation} \label{ff6}
	C_r(j,k)=\left\lbrace \begin{array}{ll}
		r_j+p_{j,k}; &   j=1,\,k=1;\\
		%{\cal R}(j,k-1)+p_{j,k}; & j=1,\,1<k\le m;\\
		\max \{r_j,\circ{\cal R}(j,k)\}+p_{j,k}; & j=1,\,1<k\le m;\\
		\circ {\cal R}(j,k)+p_{j,k}; & 1\le j\le n,\,1\le k\le m.
	\end{array}\right.
\end{equation}
\begin{equation} 	%\label {ff7}
C_{rep}(j,k)=\left \{\begin{array}{ll}
p_{j,k}; &j=1,\,k=1;\\
%{\cal R}(j-1,k)+p_{j,k}; &1< j \le n,\,k=1,\,j\,mod\,q\ne 0;\\
%{\cal R}(j-1,k)+p_{j,k}+\stackrel{rep}{\Delta}_k; &1< j \le n,\,k=1,\,j\,mod\,q=0;\\
%{\cal R}(j,k-1)+p_{j,k}; &j=1,\,1< k\le m;\\	
\circ {\cal R}(j,k)+p_{j,k}; &1\le  j \le n,\,1\le k\le m,\,j\,mod\,q\ne 0;\\
\circ {\cal R}(j,k)+p_{j,k}+\stackrel{rep}{\Delta}_k; &1\le j \le n,\,1\le k\le m,\,j\,mod\,q=0.\\
\end{array}\right.
\end{equation}

\begin{equation} 	%\label {in}
C_{in}(j,k)=\left \{\begin{array}{ll}
\stackrel{in}{\Delta}_{j,k}+p_{j,k}; &j=1,\,k=1;\\
\stackrel{in}{\Delta}_{j,k}+{\cal R}(j,k-1)+p_{j,k}; &j=1,\,1< k\le m;\\
%{\cal R}(j-1,k)+p_{j,k}; &	1<j\le n,\,k=1;\\
\circ {\cal R}(j,k)+p_{j,k}; &1< j\le n,\,1\le k\le m;\\
\end{array}\right.
\end{equation}
\begin{equation} \label{ff9}
C_{SDST}(j,k)=\left\lbrace \begin{array}{ll}
p_{j,k}+\tau_{Q(j),Q(j),k}; & j=1,\, k=1; \\
%{\cal R}(j,k-1)+p_{j,k}+\tau_{Q(j),Q(j),k}; & j=1,\, 1<k\le m;\\
%{\cal R}(j-1,k)+p_{j,k}+\tau_{Q(j-1),Q(j),k}; & 1<j \le n,\, k=1;\\
\circ {\cal R}(j,k)+p_{j,k}+\tau_{Q(j-1),Q(j),k}; &\, 1\le j \le n,\, 1\le k\le m.
\end{array}\right.
\end{equation}

For the permutation variant, the set of functions is as follows:
%$${\hat{\cal R}}\{n,m,{\bar \pi}\}:{\cal J}\times {\cal M}\to N_0\cup \{\sharp\}.$$
\begin{equation}\label{fp1}
	{\hat {\cal R}}(\alpha,k)={\hat R}_D(\alpha,k).\\
\end{equation}
%--------------------------------------------------------
%$${\hat R}_D\{n,m,\bar{D},{\bar \pi}\}:{\cal J}\times {\cal M}\to N_0\cup %\{\sharp\}.$$
\begin{align} 	\label {Dj}
	{\hat R}_{D}(\alpha,k)=\left \{\begin{array}{ll}
		{\hat R}_p(\alpha,k); &  {\hat R}_p(\alpha,k)\le D_{\pi(\alpha)},\,1\le \alpha\le n,\,1\le  k\le m;\\			
		\sharp; & {\hat R}_p(\alpha,k)> D_{\pi(\alpha)},\,1\le  \alpha\le n,\,1\le  k\le m.\\
	\end{array}\right.
\end{align}
%--------------------------------------------------------
%$${\hat R}_p\{n,m,T_{e},\stackrel{p}{\Delta},P(n,m),{\bar \pi}\}:{\cal J}\times %{\cal M}\to  N_0\cup \{\sharp\},$$
$$
%\begin{equation*} \label{f4ap}
{\hat R}_p(\alpha,k)=\left\lbrace \begin{array}{lll}
{\hat R}(\alpha,k); & {\hat R}(\alpha,k)\le T_e,&1\le \alpha\le n,\, 1\le k\le m;\\
%{\hat R}(\alpha,k); & {\hat {\cal R}}(\alpha,k-1)> T_e, & \alpha=1,\, 1<k\le m;\\
%{\hat R}(\alpha,k); & {\hat {\cal R}}(\alpha-1,k)> T_e, & 1<\alpha\le n,\, k=1;\\
{\hat R}(\alpha,k); & \circ {\hat{\cal  R}}(\alpha,k)> T_e,&1\le \alpha\le n, 1\le k\le m;\\
		% v1 -- j=1,  1<k
		%T_{e}+\stackrel{p}{\Delta}+ &{\hat {\cal  R}}(\alpha,k-1)\le T_e, &\\
	%\;\;+({\hat R}(\alpha,k)-{\hat {\cal  R}}(\alpha,k-1)); &
    %{\hat R}(\alpha,k)>T_e, &j=1,\, 1< k\le m;\\
    % v2   --  1<j,  k=1
		%T_{e}+\stackrel{p}{\Delta}+ &{\cal R}(j-1,k)\le T_e, &\\
	%\;\;+({\hat R}(\alpha,k)-{\cal R}(j-1,k)); &
   %{\hat R}(\alpha,k)>T_e, &1<j\le n,\,k=1;\\
    %  v3   -- 1<j, 1<k
T_{e}+\stackrel{p}{\Delta}+&\circ {\hat{\cal  R}}(\alpha,k)\le T_e,\\
\;\;+({\hat R}(\alpha,k)-\circ {\hat{\cal  R}}(\alpha,k)); &\,{\hat R}(\alpha,k)>T_e,&1\le \alpha\le n,\, 1\le k\le m.\\
%&\qquad\qquad\qquad\qquad  \qquad\qquad    1\le \alpha\le n,\, 1\le k\le m.\\
\end{array}\right.
%\end{equation*}
$$
%--------------------------------------------------------
\begin{equation} \label{fp4}
	{\hat R}(\alpha,k)=\left\lbrace \begin{array}{ll}
		{\hat C}(\alpha,k); & type(k)={\bf C};\\
		{\hat C}_r(\alpha,k); &  type(k)={\bf C_r};\\
		%\qquad\dots&\\
		%{\hat C}_D(\alpha,k); &  type(k)={\bf C_D};\\
		{\hat C}_{rep}(\alpha,k); & type(k)={\bf C_{rep}};\\
		{\hat C}_{in}(\alpha,k); & type(k)={\bf C_{in}};\\
		{\hat C}_{SDST}(\alpha,k); & type(k)={\bf C_{SDST}}.\\
		%{\hat C}_{p}(\alpha,k); & type(k)={\bf C_p}.
	\end{array}\right.
\end{equation}

%--------------------------------------------------------
%$$\hat{C}\{n,m,P(n,m),{\bar\pi} \}:{\cal \hat{J}}\times {\cal M}\to N_0 %\cup\{\sharp\},$$
 \begin{equation} %\label{f4p}
 {\hat C}(\alpha,k)=\left\lbrace \begin{array}{ll}
 p_{\pi(\alpha),k}; & \alpha=1,\, k=1; \\
%{\hat {\cal R}}(\alpha,k-1)+p_{\pi(\alpha),k}; & \alpha=1,\, 1<k\le m;\\
%{\hat {\cal R}}(\alpha-1,k)+p_{\pi(\alpha),k}; & 1<\alpha\le n,\, k=1;\\
 \circ {\hat {\cal R}}(\alpha,k)+p_{\pi(\alpha),k}; &\, 1\le\alpha\le n,\, 1\le k\le m.
 \end{array}\right.
 \end{equation}

%---------------------------------------------------------
%$${\hat C}_r\{n,m,\bar{r},P(n,m),\bar \pi\}:{\cal J}\times {\cal M}\to N_0\cup %\{\sharp\}.$$
\begin{equation}%\label{f3p}
	{\hat C}_r(\alpha,k)=\left\lbrace \begin{array}{ll}
		r_{\pi(\alpha)}+p_{\pi(\alpha),k}; &   \alpha=1,\,k=1;\\
		%{\hat {\cal R}}(\alpha,k-1)+p_{\pi(\alpha),k}; & %\alpha=1,\,1<k\le m;\\
		\max \{r_{\pi(\alpha)},{\hat {\cal R}}(\alpha-1,k)\}+p_{\pi(\alpha),k}; & \alpha=1,\,1<k\le m;\\
		\circ {\hat {\cal R}}(\alpha,k)+p_{\pi(\alpha),k}; & 1< \alpha\le n,\,1<k\le m.
	\end{array}\right.
\end{equation}

%----------------------------------------------------------
%$${\hat C}_{rep}\{n,m,q,\stackrel{rep}{\Delta},P(n,m),{\bar \pi}\}:{\cal J}\times %{\cal M}\to  N_0\cup\{\sharp\}.$$
\begin{equation} 	%\label {repp}
	{\hat C}_{rep}(\alpha,k)=\left \{\begin{array}{ll}
		p_{\pi(\alpha),k}; &\alpha=1,\,k=1;\\
		%{\hat {\cal R}}(\alpha-1,k)+p_{\pi(\alpha),k}; &1< \alpha \le n,\,k=1,\,
        %\alpha\,mod\,q\ne 0;\\
		%{\hat {\cal R}}(\alpha-1,k)+p_{\pi(\alpha),k}+\stackrel{rep}{\Delta}_k; &1< \alpha %\le n,\,k=1,\,\alpha\,mod\,q=0;\\
		%{\hat {\cal R}}(\alpha,k-1)+p_{\pi(\alpha),k}; &\alpha=1,\,1< k\le m;\\	
		\circ {\hat {\cal R}}(\alpha,k)+p_{\pi(\alpha),k}; &1\le \alpha \le n,\,1\le k\le m,\,
        \alpha\,mod\,q\ne 0;\\
		\circ {\hat {\cal R}}(\alpha,k)+p_{\pi(\alpha),k}+\stackrel{rep}{\Delta}_k; &1\le \alpha \le n,\,1\le k\le m,\,
        \alpha\,mod\,q=0.\\
	\end{array}\right.
\end{equation}

%----------------------------------------------------------
%$${\hat C}_{in}\{n,m,P(n,m),\|\stackrel{in}{\Delta}_{j,k}\|,{\bar \pi} \}:{\cal %J}\times {\cal M}\to N_0\cup\{\sharp\}, $$
\begin{equation} 	%\label {inp}
{\hat C}_{in}(\alpha,k)=\left \{\begin{array}{ll}
\stackrel{in}{\Delta}_{\pi(\alpha),k}+p_{\pi(\alpha),k}; &\alpha=1,\,k=1;\\
\stackrel{in}{\Delta}_{\pi(\alpha),k}+{\hat {\cal R}}(\alpha,k-1)+p_{\pi(\alpha),k}; &\alpha=1,\,1< k\le m;\\
%{\hat {\cal R}}(\alpha-1,k)+p_{\pi(\alpha),k}; &	1<\alpha\le n,\,k=1;\\
\circ {\hat {\cal R}}(\alpha,k)+p_{\pi(\alpha),k}; &1<\alpha\le n,\,1\le k\le m.\\
\end{array}\right.
\end{equation}
%--------------------------------------------------------------------------------------
%$${\hat C}_{SDST}\{n,m,z,P(n,m),{\bar \pi},\|\tau_{v,w,k}\|\}:{\cal J}\times %{\cal M}\to N_0\cup  \{\sharp\}.$$
\begin{equation} \label{sdst}
{\hat C}_{SDST}(\alpha,k)=
        \left\lbrace \begin{array}{ll}
p_{\pi(\alpha),k}+\tau_{Q(\pi(\alpha)),Q(\pi(\alpha)),k}; & \alpha=1,\, k=1; \\
%{\cal {\hat R}}(\alpha,k-1)+p_{\pi(\alpha),k}+\tau_{Q(\pi(\alpha)),Q(\pi(\alpha)),k}; & \alpha=1,\, 1<k\le m;\\
%{\cal {\hat R}}(\alpha-1,k)+p_{\pi(\alpha),k}+\tau_{Q(\pi(\alpha-1)),Q(\pi(\alpha)),k}; & 1<\alpha\le n,\, k=1;\\
\circ {\cal {\hat R}}(\alpha,k)+p_{\pi(\alpha),k}+\tau_{Q(\pi(\alpha-1)),Q(\pi(\alpha)),k}; &\, 1\le \alpha\le n,\, 1\le k\le m.
\end{array}\right.
\end{equation}
In all cases, a fixed order of jobs was considered to simplify the problem and make it more clear.
Permutation functions are of practical interest.
In what follows, only permutation functions will be considered, and it will be assumed that functions with a fixed order of jobs are a special case.
In conclusion of this section, let us return to the issue of function classification.
From the formulas (\ref{fp1}, \ref{fp4}), it is clear that they can be called ''interface'', since they do not directly affect the calculated schedule.
The remaining functions can be divided into three sets:
\begin{enumerate}
	\item 
	directly calculating the schedule\\ 
	${\hat{\cal S}}_x =\{{\hat C}(j,k),\,{\hat C}_r(j,k),\,{\hat C}_{rep}(j,k),\, {\hat C}_{in}(j,k),\, {\hat C}_{SDST}(j,k)\}$;\\
	A function belonging to a given set will be denoted by the variable 
		${\hat C}_x$, if ${\hat C}_x\in {\hat{\cal S}}_x$.
	\item 
schedule adjustments based on the constraints 
${\hat{\cal S}}_m =\{{\hat C}_p(j,k)\}$.\\
	A function belonging to a given set will be denoted by the variable 
		${\hat C}_m$, if ${\hat C}_m\in {\hat{\cal S}}_m$.
\item 
checking the schedule feasibility based on the constraints\\  
%${\cal S}_c=\{R_D(j,k)\}$ and 
${\hat{\cal S}}_c =\{{\hat R}_D(j,k)\}$.\\
	A function belonging to a given set will be denoted by a variable 
		${\hat C}_c$, if ${\hat C}_c\in {\hat{\cal S}}_c$. The union of the sets will be denoted by ${\hat{\cal C}}_g$. The function belonging to the union of the sets will be denoted by 
	%$C_g$ ò.å.  $C_g \in {\cal C}_x \cup {\cal C}_m \cup {\cal C}_c$ or 
	${\hat C}_g$, if ${\hat C}_g\in \{{\hat {\cal S}}_x \cup {\hat {\cal S}}_m \cup {\hat {\cal S}}_c$\}.
\end{enumerate}	
In this case, the second and third sets contain one function each, but there may be several functions in the set.
On the other hand, either or both of these sets may be empty.
The correct sequence of application (superposition) of these functions will be considered to be the one described above.
Within a set, functions can be applied in any order necessary to solve the problem.

Functions that are classified as interface functions, when functions are added to one of the sets or removed from it, are adjusted to include or exclude the function. However, they remain interface functions. 
\section{Objective functions}
Let us consider the description of objective functions when using recursive functions to find a schedule.
The base value is the completion time of the $j$-th job by the $k$-th machine: $C(j,k)$.
Let us consider objective functions other than the  {\em makespan}.
Let $d_j$ be the due date for completing the service.
This parameter determines the point in time by which it is desirable to complete a job. 
It is necessary to distinguish between the desired and the ultimate times of job completion (due date $d_j$ and deadline $D_j$).
The due date can be violated, although this imposes a penalty that affects the value of the objective function of the problem.

The following optimality criteria \cite{L2} and their expression through recursive functions can be formulated:
\begin{itemize}
	\item 
$L_j$ is the lateness equal to $t_j-d_j$, where $t_j$ is the moment of completing job $j$;\\
 $L_{max}\to \min$ is the criterion for minimizing the maximum lateness $$L_{max} =  \max_{1\le j\le n}\{L_j\}.$$
 In our case
 $$L_{max} =  \max_{1\le k\le m}\max_{1\le j\le n}\{C(j,k)-d_j\}.$$
 \item 
 $T_j$ is the delay equal to $\max\{0, t_j-d_j\}$ known as tardiness;\\
$\sum_{j=1}^{n}{T_j}\to \min$ is the criterion for minimizing total tardiness of the jobs  or
 $$T_{max} =  \max_{1\le k\le m} \max_{1\le j\le n}\{0,C(j,k)-d_j\}$$
corresponding to maximum tardiness.
\item 
$U_j$ -- accounting for the number of delays in the job completion, $U_j$ equals to $0$ if $t_j\le d_j$, and $1$, otherwise.
$\sum_{j=1}^{n}{U_j}\to \min$  is the criterion for minimizing the number of tardy jobs or
 $$U_{min} =  \sum_{ j=1}^{n}\{C(j,m)>d_j\}.$$
\end{itemize}

If a weight (significance) $w_j$ of job $j$ is given, the above criteria are calculated by multiplying the original value by the coefficient $w_j$.
For example, the weighted tardiness $w_j\cdot T_j$ is calculated as $w_j\cdot \max\{0, C(j,m) - d_j\}$.
\section{Recursive functions}
The main goal of this section is to characterize the set and basic properties of recursive functions that will allow: 
\begin{itemize}
	\item 
to describe the extensions of the PFSP that arise from practice;
\item 
to not go beyond the FSP complexity class;
\item 
to use a significant number of existing optimization methods.
\end{itemize}
It is known that the vast majority of PFSPs considered are NP-hard, and many solvers are used for them.
It would be wrong to create a set of recursive functions that could transfer the problem to a higher level of complexity, since in this case the complexity of the optimization algorithm increases due to the complexity of the algorithm for describing the problem and constructing the recursive functions.
In this section, we will consider the definitions and properties of recursive functions \cite{Ro,HS,St} in the part that can have an application in scheduling theory.

A function $f:N_0\to N_0$, with natural arguments and values, is called computable if there exists an algorithm $A$ that can be used to calculate the value of the function $f(n)$, that is, an algorithm $A$ such that
\begin{itemize}
	\item 
	if $f(n)$ is defined for some natural number $n$, then algorithm $A$ stops at input $n$ and prints $f(n)$;
	\item 
	if $f(n)$ is not defined, then algorithm $A$ does not stop at input $n$.
\end{itemize}
The inputs and outputs of algorithms can be not only natural numbers, but also binary strings, pairs of natural numbers, finite sequences of words, and generally any ''constructive objects''.
The concept of a computable function with two natural arguments can be defined in a similar way.

{\bf Primitive recursive functions}.
In the existing literature, the theory of recursive functions is mainly focused on problems of fundamental mathematics.
For this purpose, basis functions and operators are defined in the theory.
Operators, by definition, create new functions from existing ones. 
A function $f(x_1, . . . , x_n)$ is called {\em primitive recursive}, if it can be obtained from basis functions by a finite number of applications of the superposition and primitive recursion operators. 
Some works, such as \cite{St}, consider a simulation of logical values, arithmetic operations, accumulators, and data structures.
The recursive functions proposed for consideration are reminiscent of functional programming, when some functions are defined through others.
We will consider functions with natural arguments and values.
Functions may not be defined everywhere.

In terms of imperative programming, primitive recursive functions correspond to program blocks that use only arithmetic operations, as well as the conditional operator {\bf if-then-else} and the arithmetic loop operator {\bf for} in which the number of iterations is known at the start of the loop. 
A primitive recursive function has no reason to ''hang'' because its construction uses operators that create functions that are defined everywhere.
If the programmer begins to use the {\bf while} loop operator, in which the number of iterations is not known in advance and, in principle, can be infinite, or {\bf go to}, or changes in the {\bf for}-loop parameter inside the loop, then they move into the class of partial recursive functions.
Partial recursive functions may not be defined for some argument values.
From the point of view of imperative programming, the result of a partial recursive function can be not only a number, but also an exception or an infinite loop corresponding to an undefined value. 
Further, to describe the problems of scheduling theory, it is proposed to use recursive functions belonging to the class of primitive recursive functions.
Let us move now to the properties of functions from the set 
${\hat{\cal C}}_g$.
\subsection{Properties of recursive functions}  % ${\hat{\cal C}}_g$}
This section describes the properties of recursive functions applicable to solving flow shop problems.
The presence of these properties in a function is not a sufficient condition for the correctness of this function.
Correctness is determined taking into account the semantics of the problem.
However, if a function does not have any of the listed properties, it is incorrect, i.e., not applicable for solving flow shop problems.
The possession of these properties is a necessary but not sufficient condition.
\begin{enumerate}
	\item 
	{\bf Recursive functions ${\hat{\cal C}}_g$ are primitive recursive functions}. This statement is a limitation that is structurally imposed on the function set. 
	\item {\bf Domain area and set of values}.
 The function must be defined on the Cartesian product of the sets ${\cal J}\times {\cal M}$ and have values in the set of positive integers, i.e. 
 \begin{equation}%\label {f6}
C_g\{{\cal   P}_g\}:\;{\cal   J}\times {\cal   M} \to N_0\cup \{\sharp\}.
 \end{equation}
	-- the first argument $j\in {\cal J}$ defines the job number,\\
	-- the second argument $k\in {\cal M}$ defines the machine number,\\
	-- ${\cal P}_g$ is a set of parameters.
	\item 
	{\bf  Monotony}. 
	The function $C_g(j,k)$ increases monotonically in both arguments, i.e.
	\begin{eqnarray*}
		C_g(j,k-1) \le C_g(j,k); && 1\le j \le n,\, 1< k \le m,\\
		C_g(j-1,k) \le C_g(j,k); && 1< j \le n,\, 1\le k \le m,\\
		C_g(j-1,k-1) \le C_g(j,k); && 1< j \le n,\, 1< k \le m.
	\end{eqnarray*}

	\item 
	The following property will be formulated as a statement.  \\
{\bf Statement~1}.
%\begin{Statement}\label{u1}
	Let the values $n,\,m,\,{\cal P}$ be given for the PFSP and some order of job execution be specified.
	In this case, for any function $C_g(j,k)\, 1\le j \le n,\, 1< k \le m,$ the following is true: for any of the $n!$ fixed orders of jobs, there is no feasible schedule whose execution time $T(j,k)$ is such that $T(j,k)<C_g(j,k)$.
%\end{Statement}

	In other words, the functions $C_g(j,k)$ must be constructed in such a way that for any fixed order of jobs, they compute the minimum value.
  In \cite{KLRW}, this statement is proved for the function $C(j,k)$. 
  By the method of mathematical induction, it can be proved that the statement~1  
  %\ref{u1} 
  is also true for the remaining functions of this paper.
 Each newly defined function must also be checked for the validity of the statement~1. % \ref{u1}.
\end{enumerate}

The presence of these properties allows developing algorithms for the optimization or constraint satisfaction.

\section{Branch and bound method}
In this section, as an example, we will consider the application of the branch and bound method to a PFSP with the set of functions described above.
In order to be able to calculate a lower bound $LB$, the function $C_g(j,k)$ must satisfy the following requirement: the minimum difference in the values	
\begin{eqnarray}%\label{f7}
C_g(j,k)-C_g(j,k-1),\; \text{ for }1\le j \le n,\, 1< k \le m,\\
\label{f8}
C_g(j,k)-C_g(j-1,k)\; \text{ for }1< j \le n,\, 1\le k \le m
\end{eqnarray}
must be computable without scheduling and depends only on the values of $j$ and $k$.
In other words, it is possible to construct a function 
$$\rho: {\cal T}\times {\cal M} \to N_0$$
such that
\begin{eqnarray*}
C_g(j,k)-C_g(j,k-1)\ge\rho(j,k),\; 1\le j \le n,\, 1< k \le m,\\
C_g(j,k)-C_g(j-1,k)\ge\rho(j,k)\; 1< j \le n,\, 1\le k \le m.
\end{eqnarray*}
In addition, it is impossible to construct a function
$$h: {\cal T}\times {\cal M} \to N_0$$
such that
\begin{eqnarray*}
C_g(j,k)-C_g(j,k-1)\ge h(j,k),\; 1\le j \le n,\, 1< k \le m,\\
C_g(j,k)-C_g(j-1,k)\ge h(j,k)\; 1< j \le n,\, 1\le k \le m.
\end{eqnarray*} 
and $h(j,k)> \rho(j,k)$. 
In this case, the function $\rho(j,k)$ is a lower bound on the function $C_g(j,k)$.
\begin{table}
	\caption{Minimum function differences $C_g(j,k)$}
	\begin{center}
		\begin{tabular}{|c|c|c|c|c|}
			\hline 
			Function $C_g$ &\textnumero \; of formula& range of $j$& range of $k$&  $\rho(j,k)$ value\\ 
			\hline 
						1&2&3&4&5\\
			\hline
			$C(j,k)$&\ref{f3}& $1\le j\le n$ &$1\le k \le m$  &$p_{j,k}$\\ 
			\hline 
			$C_r(j,k)$&\ref{f3c}& $1\le j\le n$ &$1\le k \le m$  &$p_{j,k}$\\ 
			\hline 
			$C_D(j,k)$&\ref{Dj}& $1\le j\le n$ &$1\le k \le m$  &$0$\\ 
			\hline 
			%	\ref{rD}&$C_{rD}(j,k)$& $1\le j\le n$ &$1\le k \le m$  &$p_{j,k}$\\ 
			%	\hline 
		$C_{rep}(j,k)$  &\ref{rep}&$j\, mod\, q\ne 0$&$1\le k \le m$  &$p_{j,k}$\\ 
			%\hline 
			&  &$j\, mod\, q= 0$  &$1\le k \le m$  & $p_{j,k}+\stackrel{rep}{\Delta}_k$ \\  
			\hline 
			$C_{in}(j,k)$&\ref{in}&$j=1$&$1\le k \le m$  &$p_{j,k}+\stackrel{in}{\Delta}_{j,k}$\\ 
			%\hline 
			&  &$1<j$&$1\le k \le m$  &$p_{j,k}$\\  
			%	\hline 
			%\ref{buf} &$C_{buf}(j,k)$&$q+1\le j$&$k=b$&$0$\\
			%	\hline
			%	 & &$q+1> j$&$k\ne b$&$p_{j,k}$\\
			\hline 
		$C_{SDST}(j,k)$ &	\ref{sdstj}&$1\le j\le n$ &$1\le k \le m$  &$p_{j,k}+\min_{j< i\le n}\tau_{Q(j),Q(i),k}$ \\ 
			\hline 
				$C_{p}(j,k)$  &\ref{f4a}& $1\le j\le n$ &$1\le k \le m$  &$p_{j,k}$\\ 
			\hline 
		\end{tabular} 
		\label{t1}
	\end{center}
\end{table}
To demonstrate how the function $\rho$ is constructed, Table \ref{t1} is presented.
The table columns have the following meaning:
\begin{enumerate}
	\item 
	function type;
	\item 
	formula number defining the function $C_g$;
	\item 
	 range of $j$ for which $\rho(j,k)$ is valid;
	\item 
	 range of $k$ for which $\rho(j,k)$ is valid;
	\item 
value of the function $\rho(j,k)$.
\end{enumerate}
The fact that these relations are valid for the functions (\ref{f3a})-(\ref{sdst}) follows from the analysis of their definitions. 
If a definition of a new function $C_g$ is introduced and the branch and bound method is used, then it is necessary that the new function has this property.

In the case of job permutations, the function $\rho(j,k)$ will depend not on the job number, but on the job ordinal number $\alpha$, i.e., $\rho(\pi(\alpha),k)$ will be calculated.
In this case, the lower bound can be calculated using the minimum values of the function given in Table \ref{t1}.
\begin{eqnarray}\label{f14}
LB(\alpha,k)=\hat C_g(\alpha,k)+\sum_{i=k+1}^{m}{\rho(\pi(\alpha),i)}+\sum_{i=\alpha+1}^{n}{\rho(\pi(\alpha),i)}.
\end{eqnarray} 
The branching algorithm fully corresponds to the one proposed in \cite{KLRW}.
If the lower bound is calculated only for the last machine, then the following formula is used:
\begin{eqnarray}\label{f14a}
	LB(\alpha,m)=\hat C_g(\alpha,m)+\sum_{i=\alpha+1}^{n}{\rho(\pi(\alpha),m)}.
\end{eqnarray} 
The considered method of calculating the lower bounds is not the best.
Here, the solution is of a demonstration nature and illustrates the applicability of a method based on recursive functions.

{\bf Example 6}. 
%\begin{Example}
Let us consider an example of a conveyor for bottling beverages and then packaging it into boxes.
Different beverages are dispensed in different bottles.
When switching from one type of beverage to another one, it is necessary to reconfigure the equipment. 
The problem requires computing an optimal schedule of the extended PFSP for a sequence of $m=5$ machines and $n=12$ jobs.
\begin{table}
	\caption{Descriptions of beverage dispensing machines}
	\begin{center}
		\begin{tabular}{|c|l|c|l|}
			\hline 
			 Machine& &Machine &\\ 
			 number&Operation & type&Function\\ 			
%			\hline 
%			1&2&3&4&5\\
			\hline
1&Bottle feed with initial setting& $\bf {C_{in}}$ &$C_{in}(j,k)$\\ 
& for the bottle type & &\\ 
			\hline 
2&Pouring a beverage& $\bf{C}$ &$C(j,k)$\\ 
			\hline 
3&Plugging the cork& $\bf{C_{SDST}}$ &$C_{SDST}(j,k)$\\ 
			\hline 
4&Sticking the label& $\bf{C_{SDST}}$ &$C_{SDST}(j,k)$\\ 
			\hline 
5&Placing a bottle into a box and& $\bf{C_{rep}}$ &$C_{rep}(j,k)$\\ 
&  feeding the full box onto a pallet&  &\\ 
			\hline 
		\end{tabular} 
		\label{t2}
	\end{center}
\end{table}
The types of machines are as follows:\\
The 1st machine dispenses a bottle with an initial setting for the type of bottle and switches to the type of beverage.\\
The 2nd machine is pouring the beverage.\\ 
The 3rd machine corks the bottles. It needs to be adjusted to suit the bottle and cork type.\\
The 4th machine applies the label. It must be adjusted suit to the beverage and bottle type.\\
The 5th machine places bottles into boxes. When the box is full, it is sealed and placed on a pallet.

\begin{figure}[!ht]
	\centering
	\includegraphics [scale=0.6]{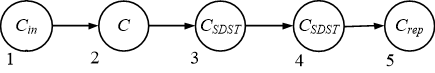}
	\caption{Example of an extended PFSP for beverage bottling.}
	\label{fig16}
\end{figure}

Figure  \ref{fig16} shows the corresponding chain of machines.
The number of production tools $z$ is equal to the number of bottle types, in this example $z=3$.
Replaceable tooling is used by machines 3 and 4, i.e. $m'=2$.
Therefore, the three-dimensional array $\tau$ has a size of $3\times 3\times 2$.
The function $Q:{\cal J}\to Z$ is defined as $Q(j)=(j-1)/m'+1$, where $x/y$ is the division of integers.
In this case, $j\in\{1..12\}$ and $Q(j)\in \{1,2,3\}$.
Permutations of jobs do not affect the value of the function.
Bottles for different types of drinks have different volumes, so the job processing times vary.
The solution to the problem is to calculate the optimal bottling schedule for 4 bottles of each of the 3 types of beverages.
The matrix of work processing times $P$, the matrix of adjustment times $\stackrel{in}{\Delta}$ and the array of tool change times $\tau$ are as follows:
$$
%	\begin{equation}
	P=\left( \begin{array}{c|ccccc}
	%	C_r&C_{in}&C_{rep}&C&C&C\\
           &1&2&3&4&5\\	
           \hline	
		1&1&3&1&1&2 \\
		2&1&3&1&1&2 \\
		3&1&3&1&1&2 \\
		4&1&3&1&1&2 \\		
		5&1&5&2&2&2\\
		6&1&5&2&2&2\\
		7&1&5&2&2&2\\
		8&1&5&2&2&2\\
		9&1&2&1&1&2 \\
		10&1&2&1&1&2 \\
		11&1&2&1&1&2 \\
		12&1&2&1&1&2 \\	
		\end{array} \right),\,
		\stackrel{in}{\Delta}=\left( \begin{array}{ccccc}
           1&2&3&4&5\\	
\hline	
		2&0&0&0&0 \\
2&0&0&0&0 \\
2&0&0&0&0 \\
2&0&0&0&0 \\	
3&0&0&0&0 \\
3&0&0&0&0 \\
3&0&0&0&0 \\
3&0&0&0&0 \\
1&0&0&0&0 \\
1&0&0&0&0 \\
1&0&0&0&0 \\
1&0&0&0&0 \\	
\end{array} \right),\,
\tau=\left( \begin{array}{ccccc}
%1&2&3&4&5\\	
%\hline	
0&5&3\\
5&0&3\\
3&5&0\\
&0&3&2\\	
&3&0&2 \\
&3&2&0 \\
\end{array} \right).		
%\end{equation}
$$
				
For machine 1, it is necessary to form a matrix of initial setup times for different jobs.
Since there are 3 batches of bottles in the example, the initial setup matrix has three time values $\{2,3,1\}$,  a different value for each batch.
On the other hand, the initial setup only applies to the first machine.
Therefore, we will describe not a matrix, but a vector of setup times.
Let us consider machine 5, which places filled bottles into a box.
The box holds 4 bottles and when the box is full, the machine packs it and places it on a pallet.
Therefore, $q=4$ and it is necessary to determine the value of the vector $\stackrel{rep}{\Delta}$.
Since there is only one machine of type $C_{rep}$ in the example, the vector will have the value $\stackrel{rep}{\Delta}=(0,0,0,0,5)$.
Here the value 5 determines the time it takes to pack the box and place it on the pallet.
The time it takes to place a bottle in a box for all types of bottles is 2. This follows from the 5th column of the matrix $P$.

\begin{figure}[!ht]
	\centering
	\includegraphics [scale=0.4]{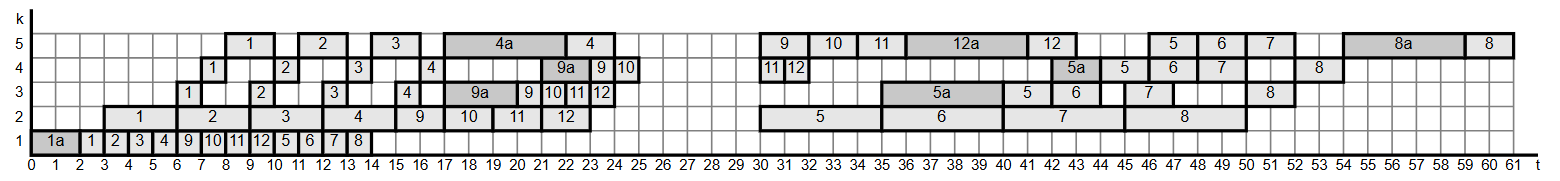}
	\caption{Time chart of the optimal beverage bottling schedule.}
	\label{fig18}
\end{figure}

Let us introduce an additional restriction: the work is interrupted by a technological break lasting $\stackrel{p}{\Delta}=5$, starting at $T_e = 25$.
The time diagram of one of the optimal schedules is shown in Figure \ref{fig18}.
The diagram shows time on the horizontal axis and machine numbers on the
vertical axis. Each rectangle in the diagram is labeled with the job
number and shows the start and end times of the job on the $k$-th
machine. Adding the letter $a$ to the job number denotes the additional
execution time $\Delta$ corresponding to the given machine.
The processing time 1a on machine 1 is determined by the initial setup.
The processing times 9a and 5a on machines 3 and 4 are determined by the transition time to a new cork and a new label. 
The processing times 4a, 12a, and 8a on machine 5 are determined by the timing of periodic operations with the box.
Of the many schedules, those in which bottling was not carried out in batches were discarded, otherwise bottles with different beverages would have ended up in the same box.
%\end{Example}

\section{Conclusion}
Since a recursive function is, by definition, calculated using some algorithm, the paper uses several examples to show a method for localizing the specifics of a problem in the body of a function. 
It should be noted that the considered choice of PFSP extensions is due to two main reasons.
First, all the described extensions are directly related to real problems of scheduling theory, i.e., they are relevant.
Second, the described extensions have a simple formalization, which makes it possible to describe the implementation of the solution quite clearly, i.e., they are demonstrative in nature.

The article presented a three-level structure of function superpositions: functions associated with the type of machine, functions for adjusting the schedule in accordance with the constraints, and functions for checking the feasibility of the schedule in accordance with the constraints.
Each of these levels can be supplemented with appropriate functions and describe a new extension of the problem.
The function implementing the predicate $and$ in \cite{KLRW} can be considered a zero-running-time machine.

Recursive functions of this type have a compact and formal description and can themselves perform the function of classifying a set of problems, i.e., the description of the set of functions and some informational support can state the  formulation of the problem.

Each of the above extensions can be significantly complicated to further approach the practice, and at the same time they will remain within the framework of the described methodology.
For example, the described extension of the work pause, depending on the technological features, may more subtly take into account the conditions of work interruption. 
This applies to other functions as well. 

Extensions such as those discussed in the article are formally described as recursive functions and provide an opportunity for further comprehensive research with comparisons of results and can increase the practical significance of the PFSP.
 
The paper does not discuss the efficiency of calculating functions.
However, it is obvious that the implementation should involve storing the results to avoid repeated calculations.

\section*{Author contributions}
Conceptualization, B.K. and A.L.; 
methodology, B.K. and A.L.;
software, A.R.;
validation, B.K., A.L., A.R. and F.W.;
formal analysis, F.W.;
investigation, B.K.; 
writing---original draft preparation, B.K.; 
writing---review and editing, F.W.; 
visualization, B.K. and A.R.; 
supervision, B.K. and A.L.; 
project administration, A.L.; 
funding acquisition, F.W. All authors have read and agreed to the published version of the manuscript.

\section*{Conflict of interest}
The authors declare no conflicts of interest.
%%%%%%%%%%%%%%%%%%%%%%%%%%%%%%%%%%%%%%%%%%
%%%%%%%%%%%%%%%%%%%%%%%%%%%%%%%%%%
%\begin{adjustwidth}{-\extralength}{0cm}

%=====================================

%For more questions regarding reference style, please refer to the \href{http://www.ncbi.nlm.nih.gov/books/NBK7256/}{Citing Medicine}.

%\section*{Supplementary (if necessary)}
\end{document}